\RequirePackage{fix-cm}
\documentclass[smallextended]{svjour3}       
\smartqed  
\usepackage{graphicx}
\usepackage[T1]{fontenc}
\usepackage[utf8]{inputenc}
\usepackage{amsmath,amsfonts}
\usepackage{pst-all}
\usepackage{graphicx}
     \usepackage{pstricks-add}

\newcommand{\N}{\mathbb N}
\newcommand{\R}{\mathbb R}
\newcommand{\T}{\mathbb T}

\renewcommand{\P}{\mathbb P}

\newcommand{\E}{\mathbb E}

\journalname{To appear in the Journal of Statistical Physics}
\begin{document}

\title{On rumour propagation among sceptics
}

\author{Farkhondeh Alsadat~Sajadi and  Rahul~Roy  }


\institute{Farkhondeh Alsadat~Sajadi  \at
              University of Isfahan, Isfahan\\
              \email{f.sajadi@sci.ui.ac.ir}         
           \and
           Rahul~Roy  \at
              Indian Statistical Institute, Delhi\\
							\email{rahul@isid.ac.in} 
}

\date{Received: date / Accepted: date}

\maketitle

\begin{abstract}
Junior, Machado and Zuluaga (2011) studied a model to understand the spread of a rumour. Their model  consists of individuals situated at the integer points of the line $\N$. An individual at the origin $0$ starts a rumour and passes it to all individuals in the interval $[0,R_0]$, where $R_0$ is a non-negative random variable. An individual located at $i$ in this interval receives the rumour and transmits it further among individuals in $[i, i+R_i]$ where $R_0$ and $R_i$ are i.i.d. random variables. The rumour spreads in this manner.  An alternate model considers individuals seeking to find the rumour from individuals who have already heard it. For this s/he asks individuals to the left of her/him and lying in an interval of a random size.  We study these two models, when the individuals are more sceptical and they transmit or accept the rumour only if they receive it from at least two different sources. 

In stochastic geometry the equivalent of this rumour process is the study of coverage of the space $\N^d$ by random sets. Our study here extends the study of coverage of space and considers the case when each vertex of $\N^d$ is covered by at least two distinct random sets.

\end{abstract}

\section{Introduction}
Gilbert (1961) introduced a model of transmission of information. This model consisted of a signal being transmitted through a relay of transmitters to its recipient. Variants of this process were later introduced and described in Maki and Thompson (1973). Two such versions are the Poisson Boolean model and the rumour processes. We briefly describe these processes here and indicate some literature connected with them.

\vspace{.3cm}

\noindent{\sc Poisson Boolean model: } The Poisson Boolean model consists of a homogenous Poisson point process $\Xi := (\xi_1, \xi_2, \ldots )$ on $ \R^d$ of intensity $\lambda$. Along with this is an independent collection of i.i.d. positive real valued random variables $\{\rho_1, \rho_2, \ldots\}$. The covered region of the Boolean model is defined to be the random set $C:= \cup_{i=1}^\infty B(\xi_i, \rho_i)$, where $B(\xi, \rho)$ is the closed ball centred at $\xi$ and of radius $\rho$ in the Euclidean norm.

This Boolean model is used to study coverage properties of space in stochastic geometry. Matheron (1968) used this Boolean model to  study natural images and their occlusion. Later, other geometric properties were studied, see e.g.,  Hall (1988) and Chiu, Stoyan, Kendall and Mecke (2013) for further details. Kertesz and Vicsek (1982) picked up this model as a natural extension on the continuum space $\R^d$ of unoriented bond/site percolation in physics. The percolation parameter being the intensity $\lambda$ with the radius random variables $\rho_1, \rho_2, \ldots$ being either constants or of a fixed distribution. A review of the mathematical details of this percolation model and its variant, the random connection model, may be seen in Meester and Roy (1996) and Penrose (2003). Gupta and Kumar (1998) used this model to study questions of signal-to-interference-ratio (SINR) and other such problems in wireless transmission, see Franciscceti and Meester (2007) for a review. 

\vspace{.3cm}

\noindent{\sc Rumour process: } Sudbury (1985) started the mathematical study of this variant of the information-transmission model introduced in Maki and Thompson (1973). Subsequently, Junior, Machado and Zuluaga (2011) christened this model as the `firework process' and introduced a different variant the `reverse firework process'. Both these models are defined on a discrete graph as opposed to the Boolean model. For convenience we describe them for the non-negative integer line $\N$. 

\noindent {\it Firework process: } The {\it homogenous}\/ firework process consists of a sequence of non-negative integer valued i.i.d. random variables $\{R_i : i \geq 0\}$. At time $0$ the origin starts a rumour and passes it onto all individuals in the interval $[0, R_0]$, and at time $t$, all individuals who heard the rumour for the first time at time $t-1$  spread the rumour, with the individual at site $j$ spreading it among all individuals in the region $[j, j+R_j]$. Note that allowing $\P\{R_j = 0\} >0$ ensures that there are individuals who are inactive. The {\it heterogeneous}\/ firework process removes the restriction on the random variables $\{R_i : i \geq 0\}$ of being identically distributed.

\noindent {\it Reverse firework process: } The reverse firework process consists of the origin who knows the rumour at time $0$, and at time $t$ an individual located at site $j$ listens to individuals in the interval $[j-R_j, j]$. If there is an individual at a site in this interval who has heard the rumour by time $t-1$, then the individual at site $j$ gets to know the rumour. Here the random variables $\{R_i : i \geq 0\}$ are as before, and accordingly we have a homogenous or heterogeneous version of the model.

Both these models have been extensively studied as models for spread of information on networks. 
Gallo, Garcia, Junior and Rodr\'{i}guez (2014) exploit a relationship between the rumour processes
and a specific discrete time renewal process to obtain various
results. Junior, Machado and Zuluaga (2014) study a version of  the firework process on homogeneous trees.
Bertacchi and Zucca (2013) study these processes in a random environment. Junior, Machado and Ravishankar (2016) summarises the current state of the research on these topics.

\vspace{.3cm}

In this paper, we study the spread of rumours among sceptics. Sceptics being those individuals who need at least $2$ different individuals from whom they receive the rumour before transmitting/accepting. The results are valid for sceptics who need the information from $k$ or more sources for acceptance or transmission, however for the sake of simplicity we restrict ourselves to $k=2$. 

For Boolean models vis-\`{a}-vis coverage of space, Athreya, Roy and Sarkar (2004) introduced a notion of `eventual coverage' which, for 1-dimension,  is equivalent to the notion of  percolation of the firework process as will be explained in the next section. This study was not restricted to the half-line, but more generally to a quadrant in $d$-dimensional space. Our paper extends this study to `double coverage', i.e., the region in the quadrant which is covered by at least $2$ distinct random sets.

In the next section we present the formal set-up of the questions we study, and state the results. In Sections 3 and 4 we prove the propositions, and in Section 5 we state and explain in brief the version of the coverage result for the Poisson Boolean model.

\section{The models and the statement of results}
We present an alternate but equivalent formulation of the homogenous firework and the reverse firework processes.
Let $X_1, X_2, \ldots$ be i.i.d.  Bernoulli ($p$) random variables, i.e.,
\begin{equation}
\label{Bern_p}
X_{\mathbf i}= \begin{cases}
       1 & \text{ with probability } p\\
       0 & \text{ with probability } 1-p.
      \end{cases}
\end{equation}
Also let $\{\rho_i: i \geq 1\}$ be a collection of i.i.d. $\mathbb N$ valued random variables, independent of the collection $\{X_i: i \geq 1\}$. Let $\rho$ denote a generic random variable with the same distribution as $\rho_i$.
In addition, let $\rho_0$ be an independent $\mathbb N$ valued random variable, independent of the collections $\{X_i: i \geq 1\}$ and $\{\rho_i: i \geq 1\}$, with $\rho_0$  having the same distribution as $\rho$.
 Taking $X_0 \equiv 1$, consider the regions
\begin{align}
\label{r:fireC}
C&:=  \bigcup_{\{i \geq 0 : X_{i} = 1\}} [{i}, {i}+\rho_i],\nonumber\\
C_{\text{rev}} &:=  \bigcup_{\{i \geq 0 : X_{i} = 1\}} [i - \rho_i , {i}].
\end{align}
The firework processes presented in the earlier section may be seen to be equivalent to this formulation by taking 
$p= \P\{R_0 > 0\}$ and $\rho$ to have the same distribution as that of $R_0 | R_0 > 0$.

\begin{remark}
\label{r:rmk1}
A simple argument using Kolmogorov's $0$-$1$ law yields that $C$ (or  $C_{\text{rev}}$) being an unbounded connected region with positive probability is equivalent to $C$ (or  $C_{\text{rev}}$) containing a region $[t, \infty)$, for some $t \geq 1$ with probability $1$.
\end{remark}

We study the spread of the rumour among sceptics. If individual at $i$ receives the rumour from at least two distinct sources, then  s/he transmits it to all individuals in the region $[i, i+\rho_{ i}]$, i.e.,
\begin{itemize}
\item[({\bf A})] for the firework process, there are two individuals at $j$ and $k$ (say) with $j \neq k$, $X_j = X_k = 1$ and $-1 \leq j, k <  i$ such that $i \in [j, j + \rho_j] \cap [k, k+ \rho_k]$; here we assume that there is an individual at location $-1$, who spreads the rumour in the region $[-1, -1+ \rho_{-1}]$, where $\rho_0$ and $\rho_{-1}$ are two independent random variables, independent of all other random variables, and each having the same distribution and $X_0 = X_{-1} \equiv 1$; 
\item[({B})] for the reverse firework process, there are two individuals at $j$ and $k$ (say) with $j \neq k$, $X_j = X_k = 1$ and $-1 \leq j, k <  i$ such that $j,k \in [{i} - \rho_{ i}, {i}]$, where $\rho_0$ and $\rho_{-1}$ are two independent random variables, independent of all other random variables, and each having the same distribution as $\rho$ and $X_0 = X_{-1} \equiv 1$.
\end{itemize}

Towards this end we define the regions 
\begin{align}
\label{r:fireD}
D &:= \{x \in \R: \text{there exist } j, \;  k \geq -1 \text{ with } j \neq k , \;  X_j= X_k = 1\nonumber \\
& \qquad \text{ and } x \in([j, j + \rho_{ j}] \cap [k, k+ \rho_{ k}])\},\nonumber\\
D_{\text{rev}} &:= \{x \in \R_+:  x \in [{i} - \rho_{ i}, {i}] \text{ for some } i \text{ with } X_{i} = 1 \text{ and there exist } \nonumber\\
& \qquad  -1\leq j \neq  k < i \text{ with }  X_j= X_k = 1 \text{ and }  j, \;  k \in [{i} - \rho_{ i}, {i}]\}.
\end{align}
We look for conditions  on the processes $\{X_i: i \geq 1\}$ and $\{\rho_i: i \geq 1\}$ such that the regions $D$ and $D_{\text{rev}}$ are unbounded connected regions with positive probability.
\begin{remark}
\label{r:rmk2} In the set-up {\rm ({\bf A})} (or in {\rm ({\bf B})}), as in the case of rumour propagation in $C$ (or $C_{\text{rev}}$), a simple argument using Kolmogorov's $0$-$1$ law yields that $D$ (or  $D_{\text{rev}}$) being an unbounded connected region with positive probability is equivalent to $D$ (or  $D_{\text{rev}}$) containing a region $[t, \infty)$, for some $t \geq 1$ with probability $1$. Remark \ref{r:rmkequivdef} given in the next section spells out the details.
\end{remark}

Thus we define 
\begin{definition}
The firework process (respectively, reverse firework process) percolates among sceptics if,  with probability $1$,  $D$ (respectively,  $D_{\text{rev}}$) contains a region $[t, \infty)$, for some $t \geq 1$.
\end{definition}
\begin{remark}
\label{r:rmk-10}
Using this equivalent definition obtained by the tail event properties allows us to study the model without considering the influence of the two initiators $X_{-1}$ and $X_0$. Indeed their influence will only be upto $\max\{\rho_{-1}, \rho_0\}$.
\end{remark}
\begin{proposition}
\label{r:even_main1}
Let 
$$
l := \liminf_{j \to \infty} j \mathbb{P}(\rho \geq j) > 1 \text{ and } L:= \limsup_{j \to \infty} j \mathbb{P}(\rho \geq j) < \infty.
$$ 
We have that the firework process given in {\rm ({\bf A})} percolates among sceptics if $p > 1/l$
and does not percolate if $p < 1/L$, i.e., 
$$
\P\{ D \supseteq [t, \infty) \text{ for some finite } t \} = \begin{cases} 1 & \text{ if }p > 1/l\\
0 & \text{ if }p < 1/L.
\end{cases}
$$
\end{proposition}

For the reverse firework process we have
\begin{proposition}
\label{r:even_main11}
For $p > 0$, the reverse firework process given in  {\rm ({\bf B})}  percolates among sceptics if and only if 
$\E(\rho) = \infty$, and, in case percolation occurs  
$\P\{ D_{\text{rev}} = \N^d\} =1$.
\end{proposition}

We note here that these are the same conditions that Junior, Machado and Zuluaga (2011) obtain for  percolation among `non-sceptical' individuals. Indeed, the above propositions also go through among more radical sceptics, i.e. if individuals need to receive the rumour from $k \geq 1$ distinct sources before they transmit/believe  the rumour.

As noted in the review article Junior, Machado and Ravishankar (2016), the above model is related to the study of coverage processes in stochastic geometry. Athreya, Roy and Sarkar (2004) introduce a notion of `eventual coverage' which, for 1-dimension,  is identical to the equivalent formulation of the percolation of rumour process as given in Remark \ref{r:rmk1}. We state the model in brief here and present the results obtained. Their formulation considers 
$\{X_{\mathbf i}: {\mathbf i} \in \mathbb N^d\}$ to be a collection of i.i.d.  Bernoulli ($p$) random variables and $\{\rho_{\mathbf i}: {\mathbf i} \in \mathbb N^d\}$  a collection of i.i.d. $\mathbb N$ valued random variables, independent of the collection $\{X_{\mathbf i}: {\mathbf i} \in \mathbb N^d\}$. 
Let $\rho$ denote a generic random variable with the same distribution as $\rho_{\mathbf i}$ and  
$$
{\mathbf C}:= \cup_{\{{\mathbf i} : X_{\mathbf i} = 1\}} ({\mathbf i}+[0,\rho_{\mathbf i}]^d)
$$ 
denote the {\it covered region}\/ of $\mathbb N^d$; here and subsequently ${\mathbf i}+[0,\rho_{\mathbf i}]^d = [i_1, i_1 + \rho_{\mathbf i}] \times \cdots \times [i_d, i_d + \rho_{\mathbf i}]$, where ${\mathbf i} = (i_1, \ldots, i_d)$.
\begin{definition}
\label{r:evencov}
 $\mathbb N^d$ is eventually covered if there exists  $ {\mathbf t} \in \mathbb N^d$ such that ${\mathbf t} + \mathbb N^d \subseteq {\mathbf C}$. 
 \end{definition}
 
\begin{remark}
\label{r:ARSJMZ}
From Remark \ref{r:rmk1}, the above definition may be seen to be equivalent to percolation of the homogenous firework process for $d=1$, and in that sense, it extends the definition of percolation for a homogenous firework process in $\N^d$.
\end{remark}
 
For our purposes we define
\begin{definition}
\label{r:evendoubly}
 $\mathbb N^d$ is {\it eventually  doubly covered}\/ if there exists  $ {\mathbf t} \in \mathbb N^d$ such that ${\mathbf t} + \mathbb N^d \subseteq {\mathbf D}$, where
\begin{align*}
{\mathbf D} &:= \{{\mathbf x} \in \mathbb R^d: \text{there exist } {\mathbf i}, \; {\mathbf j}\in \mathbb N^d \text{ with } {\mathbf i} \neq {\mathbf j} \text{ and } X_{\mathbf i} = X_{\mathbf j} = 1 \\
& \qquad \text{ such that } x \in({\mathbf i} + [0, \rho_{\mathbf i}]^d) \cap ({\mathbf j}+[0,\rho_{\mathbf j}]^d)\}.
\end{align*}
\end{definition}
We have 
\begin{proposition}
\label{r:even_main2}
(i) For $d=1$, let 
$$
l := \liminf_{j \to \infty} j \mathbb{P}(\rho \geq j) > 1 \text{ and } L:= \limsup_{j \to \infty} j \mathbb{P}(\rho \geq j) < \infty.
$$ We have 
$$
\mathbb{P}_p (\mathbb N \text{ is eventually doubly covered}) = \begin{cases} 
1 & \text{if } p > 1/l\\
0 & \text{if } p < 1/L.
\end{cases}
$$
(ii) For $d\geq 2$ and $p > 0$, we have 
$$
\mathbb{P}_p (\mathbb N^d \text{ is eventually doubly covered}) = \begin{cases} 
1 & \text{if } \liminf_{j \to \infty}  j \mathbb{P}(\rho \geq j) > 0\\
0 & \text{if } \lim_{j \to \infty}  j \mathbb{P}(\rho \geq j) =0.
\end{cases}
$$
\end{proposition}

Apropos the reverse firework process in higher dimensions, for $\{X_{\mathbf i}: {\mathbf i}\in \N^d\}$ and $\{\rho_{\mathbf i} : {\mathbf i} \in \N^d\}$ as earlier, let
$$
{\mathbf C}_{\rm{rev}} := \cup_{\{{\mathbf i} : X_{\mathbf i} = 1\}} ({\mathbf i}+[-\rho_{\mathbf i}, 0]^d)
$$
and for $-{\mathbf 1} = (-1,-1, \ldots, -1)$, ${\mathbf 0} = (0,0,\ldots , 0)$, $X_{\mathbf 1} = X_{\mathbf 0} = 1$, 
\begin{align*}
{\mathbf D}_{\rm{rev}} := & \left\{{\mathbf x} \in \R_+^d : {\mathbf x} \in ({\mathbf i}+[-\rho_{\mathbf i}, 0]^d) \text{ for some } 
{\mathbf i} \text{ with } X_{\mathbf i} = 1\right.\\
& \text{ and there exist } {\mathbf j}\neq {\mathbf k}, \; {\mathbf j} , {\mathbf k} \in \left( \N^d \cup \{- {\mathbf 1}, {\mathbf 0}\}\setminus \{ \mathbf i\} \right) \\
& \left. \text{such that } X_{\mathbf j} = X_{\mathbf k} = 1 \text{ and } {\mathbf j}, {\mathbf k} \in ({\mathbf i}+[-\rho_{\mathbf i}, 0]^d) \right\}.
\end{align*}

\begin{proposition}
\label{r:even_main22}
For $p > 0$, the reverse firework process on $\N^d$  percolates among sceptics if and only if 
$\E(\rho^d) = \infty$, and, in case percolation occurs  
$\P\{ D_{\text{rev}} = \N^d\} =1$.
\end{proposition}
Finally, because of the binomial approximation of the Poisson process, Proposition \ref{r:even_main2} has a natural extension to Poisson processes. This is relegated to the last section of this paper.

In stochastic geometry the notion of coverage of space has received extensive
 attention. In particular
Hall (1988) and Chiu, {\it et al}\/ (2013) provide a review of the topics studied. Our endeavour in this paper may be viewed as an effort to introduce a notion of `reinforced coverage'.

\section{Proofs of Propositions \ref{r:even_main1} and \ref{r:even_main11} }
\noindent{\bf{Proof of Proposition \ref{r:even_main1}:} }
Let the random regions $C$ and $D$ be as in (\ref{r:fireC}) and (\ref{r:fireD}) respectively. 
In view of Remark \ref{r:rmk-10}, we may simplify the process to be defined only for the positive integers and ignore the individuals located at $0$ and $-1$. Indeed, for any sample point $\omega$, if $D(\omega)$ contains an interval $[t, \infty)$ for some finite $t$, it will also contain the interval $[\max\{t, \rho_0(\omega), \rho_1(\omega)\}, \infty)$, see also Remark \ref{r:rmkequivdef} given later.
Thus with a slight abuse of notation, we define $C$ and $D$ from only site $1$ onwards.

For $p > 0$ and $i \geq 1$,  let
$$
A_i := \{i \not\in C\}
 \text{ and } 
B_i := \{i \not \in D\}.
$$
Taking $G(i) = \P(\rho \geq i)$ and $g_p(i) = 1-p G(i)$, we observe that
\begin{align}
\label{PB_i}
& \mathbb{P}_p(A_i)
= \prod_{l=0}^{i-1}\P_p\left( \{X_{i-l} = 0\} \cup \{X_{i-l} = 1 \text{ and } \rho_{i-l} < l\}\right)  =\prod_{l=0}^{i-1} g_p(l), \nonumber\\
&\mathbb{P}_p(B_i) \nonumber\\
=&\mathbb{P}\left( A_i \cup \{(\text{there exists exactly one $j$ with $X_j = 1$ such that } i \in [j, j+ \rho_j] \}\right)\nonumber\\
=& \mathbb{P}_p(A_i) + \sum_{l=0}^{i-1} \mathbb{P}_p(X_{i-l} = 1, \; i \leq i-l + \rho_{i-l}, \; i \not\in \cup_{\{j \neq i- l, X_j = 1\}} [j, j+ \rho_j])\nonumber\\
= &\mathbb{P}_p(A_i) + p \prod_{l=1}^{i-1} g_p(l)+ p(1-p)\sum_{k=1}^{i-1} G(k)\prod_{l\neq k, l=1}^{i-1} g_p(l) .
\end{align}

Now suppose $p > 1/l$, where $l$ is as in Proposition \ref{r:even_main1}. We will show that 
\begin{align}
\label{B_sum_infty}
\sum_i \mathbb{P}_p(B_i) < \infty, 
\end{align}
 which, by Borel-Cantelli lemma yields  $\mathbb{P}_p(B_i \text{ occurs finitely often}) = 1$, i.e. there is a random variable $T$, with $T<\infty$ almost surely, such that $\P_p\{D \supseteq [T, \infty)\} = 1$.
 
From  the proof of Proposition 3.1 (a) of Athreya, Roy and Sarkar (2004) we have that  
\begin{align}
\label{r:cond0}
\sum_i  \prod_{l=1}^{i-1} g_p(l) < \infty \text{ for } p > 1/l.
\end{align}
 Thus, from expression (\ref{PB_i}), to show (\ref{B_sum_infty}),  we need to show
\begin{align}
 \label{r:after}
\sum_{i=1}^\infty \sum_{k=1}^{i-1} G(k)\prod_{l\neq k, l=1}^{i-1} g_p(l)  < \infty . 
\end{align}

Towards this end we note that since $p > 1/l$ there exists $\eta > 1$ such that $l > \eta /p$.
Also fix $i_0$  such that  $iG(i) > \frac{\eta}{p}$ for all $i \geq i_0$.
Now 
\begin{align*}
\sum_{k=1}^{i-1} G(k)\prod_{l\neq k, l=1}^{i-1} g_p(l) & =  
 \prod_{l=i_{0}+1}^{i-1} g_p(l)\left\{\sum_{k=1}^{i_0} G(k)\prod_{l\neq k,l=1}^{i_0} g_p(l)\right\}\\
& \quad + \prod_{ l=1}^{i_0} g_p(l) \left\{\sum_{k={i_0}+1}^{i-1} G(k)\prod_{l\neq k,l=i_{0}+1}^{i-1} g_p(l)\right\}.\\
\end{align*}
From our observation (\ref{r:cond0})  and  since $\sum_{k=1}^{i_0} G(k)\prod_{l\neq k,l=1}^{i_0} g_p(l)$ is a constant for fixed $i_0$, to show (\ref{r:after})
 we only need to show
\begin{align}
\label{B_fin}
\sum_{i= i_0 + 1}^\infty \sum_{k={i_0}+1}^{i-1} G(k)\prod_{l\neq k,l=i_{0}+1}^{i-1} g_p(l) & =\sum_{k=i_{0}+1}^{\infty} G(k)\sum_{i=k+1}^{\infty}\prod_{l\neq k,l=i_{0}+1}^{i-1} g_p(l)\nonumber\\
& < \infty .
\end{align}

In preparation for the ratio test we will use to show (\ref{B_fin}), let
\begin{align*} 
b_k &: = G(k)\sum_{i=k+1}^{\infty}\prod_{l\neq k,l=i_{0}+1}^{i-1} g_p(l)\\
& = G(k)\prod_{l=i_{0}+1}^{k-1} g_p(l)\left\{ 1+ \sum_{i=k+2}^{\infty}\prod_{l=k+1}^{i-1} g_p(l)\right\}.
\end{align*} 
Also
\begin{align}
\label{r:bratio}
\frac{b_{k+1}}{b_{k}}&=\frac{G(k+1)}{G(k)}g_p(k)\frac{1+ \displaystyle\sum_{i=k+3}^{\infty}\displaystyle\prod_{l=k+2}^{i-1} g_p(l)}{1+ \displaystyle\sum_{i=k+2}^{\infty}\displaystyle\prod_{l=k+1}^{i-1} g_p(l)}.
\end{align}

Since we will use such computations again in the next section, we elaborate the details here.
For the numerator in the above expression, first note that
\begin{equation}
\label{r:conda}
\frac{G(k+1)}{G(k)} \leq 1  \text{ and } g_p(k) \leq 1 \text{ for all } k.
\end{equation}
From (\ref{r:cond0}) we have $\sum_{i=k+3}^\infty\prod_{l=1}^{i-1} g_p(l) = \prod_{l=1}^{k+1} g_p(l)
\sum_{i=k+3}^\infty\prod_{l=k+2}^{i-1} g_p(l)
< \infty$ and thus,  $ \sum_{i=k+3}^\infty\prod_{l=k+2}^{i-1} g_p(l) < \infty$ for fixed $k \geq 1$, i.e., 
\begin{equation}
\label{r:condb}
\sum_{i=k+3}^\infty\prod_{l=k+2}^{i-1} g_p(l) < 0.1 \text{ for all } k \text{ large enough.}
\end{equation}
Hence, till now we have
\begin{align}
\label{r:bratio1}
\frac{b_{k+1}}{b_{k}} < (1.1) \left[1+ \displaystyle\sum_{i=k+2}^{\infty}\displaystyle\prod_{l=k+1}^{i-1} g_p(l)\right]^{-1} \text{ for all }  k \text{ large enough.}
\end{align}
Writing $\displaystyle\sum_{i=k+2}^{\infty}\displaystyle\prod_{l=k+1}^{i-1} g_p(l)= g_p(k+1)\left[1+\displaystyle\sum_{i=k+3}^{\infty}\displaystyle\prod_{l=k+2}^{i-1} g_p(l)\right]$,
and noting that $G(k)\downarrow 0$ as $k \uparrow \infty$, we have
\begin{equation}
\label{r:condc} 
0.9 < g_p(k+1) \leq 1 \text{ for all } k \text{ large enough.}
\end{equation}
Since $\displaystyle\sum_{i=k+3}^{\infty}\displaystyle\prod_{l=k+2}^{i-1} g_p(l) \geq 0$, we have from (\ref{r:bratio1}) and (\ref{r:condc}),
$$
\frac{b_{k+1}}{b_{k}} < \frac{1.1}{1.9} \text{ for all }  k \text{ large enough.}
$$
The ratio test  proves (\ref{B_fin})
which shows that (\ref{B_sum_infty}) holds.

From Proposition 3.1 (b) of Athreya, Roy and Sarkar (2004), we know $\mathbb{P}_p (C \supseteq [t, \infty) \text{ for some $t$ finite}) = 0$ for $p < 1/L$ and  $L:= \limsup_{j \to \infty} j \mathbb{P}(\rho > j) < \infty$. Now 
$\{D\supseteq [t, \infty) \text{ for some $t$}\} \subseteq
\{C \supseteq [t, \infty) \text{ for some $t$}\}$ which completes the proof of Proposition \ref{r:even_main1}. \hfill $\square$
\begin{remark}
\label{r:rmkequivdef}
The Borel-Cantelli argument used in the proof gives that the random variable $T$ defined at the beginning of this section is finite almost surely.
Thus, given $\epsilon >0$, we may obtain $m$ such that $\P(T \leq m) > \epsilon$. Now the probability that $X_i = 1$ for $i = 1, \ldots , m$ is $p^m > 0$. Hence, an application of the FKG inequality yields 
\begin{align*}
& \P_p\{\text{the firework process percolates among sceptics}\}\\
& \geq \P_p(\{D\supseteq [m, \infty)\} \cap \{X_1 = \cdots = X_m = 1\}) > \epsilon p^m > 0.
\end{align*}
This corroborates  Remark \ref{r:rmk2}.

\end{remark}

\vspace{.3cm}
\noindent{\bf{Proof of Proposition \ref{r:even_main11}:}}
We present an argument here which also simplifies the proof of Proposition 2.3 (i) of Junior, Machado and Zuluaga (2011).

Fix $p > 0$. Let $A_n:= \{X_n = 1 \text{ and } \rho_n \geq n+1\}$. Thus, if $A_n$ occurs, the individual at site $n$ can access the rumour from the two individuals at sites $-1$ and $0$. Noting that (i) $\{A_n: n\geq 1\}$ is a collection of independent events, (ii) $\P_p(A_n) = p \P(\rho \geq n+1)$
and (iii) $\E( \rho) = \infty$ implies that $\sum_{n\geq 1} \P(\rho \geq n+1) = \infty$, we have from Borel-Cantelli lemma, $\P_p(A_n \text{ occurs infinitely often}) = 1$. Thus  $\P_p (D_{\text{rev}} = \N) = 1$ whenever $\E( \rho) = \infty$. 

To complete the proof of Proposition \ref{r:even_main11} we note that Proposition 2.3 (ii) of Junior, Machado and Zuluaga (2011) says $\P_p(C_{\text{rev}} \supseteq [t,\infty) \text{ for some } t) = 0$ whenever 
$\E( \rho) <\infty$. Since $D_{\text{rev}} \subseteq C_{\text{rev}}$, the proof is complete. \hfill $\square$
\begin{remark}
\label{r:sceptic} The proof of Proposition \ref{r:even_main11} above reinforces our contention in Section 1 that our results are valid for any `radical' sceptic.
Indeed,  whenever $\E( \rho) = \infty$  we have $\P\{$there is an infinite increasing sequence $n_1, n_2, \ldots$ such that for any $k$, the individuals $n_j, \; j \geq k$ have heard the rumour from at least $k$ different sources$\} = 1$.  However for extending the proof of Proposition \ref{r:even_main1} a bit more work has to be done. \end{remark}



\section{Proofs of Propositions \ref{r:even_main2} and  \ref{r:even_main22}}
\noindent{{\bf Proof of Proposition \ref{r:even_main2}:} 
As mentioned earlier  Proposition \ref{r:even_main2} (i) is just Proposition \ref{r:even_main1} rephrased. Thus we need to prove Proposition \ref{r:even_main2} (ii).

First note that from Proposition 3.2(a) of Athreya, Roy and Sarkar (2004), we know that  if $\lim_{j\longrightarrow\infty} j \mathbb{P}(\rho \geq j)=0 $ then $\mathbb P_p(\mathbb N^d \text{ is eventually covered}) = 0$, and so  $$\mathbb P_p(\mathbb N^d \text{ is eventually doubly covered}) = 0.$$

We prove Proposition \ref{r:even_main2} (ii) for the case $d=2$; the proof carries through in a straightforward fashion for higher dimensions. Since the proof is technical and long, we break it up into various steps.

\vspace{.3cm}

\noindent {\sc Step 1: Prelude} 

\vspace{.3cm}

Fix $0 < p < 1$ and we assume that $\liminf_{j\longrightarrow\infty} j \mathbb{P}(\rho \geq j)>0 $.

For $(i,j) \in \mathbb N^2$, define
$$
B_{i,j} := \{(i,j) \not \in {\bf D}\}.
$$
If we show that, for some $N\geq 1$,$$
\sum_{{i,j} \geq N}\P _p(B_{i,j}) < \infty, 
$$
then Borel-Cantelli lemma guarantees that, with probability $1$, there exists $N_0 \geq 1$ such that 
$(i,j)  \in {\bf D}$ for all $i,j \geq N_0$, i.e. we have eventual coverage. Also an argument similar to that in Remark \ref{r:rmkequivdef} shows that, in this case, percolation occurs among sceptics in $d=2$.
\vspace{.3cm}
\begin{figure}[htb]
		\includegraphics[scale=0.5]{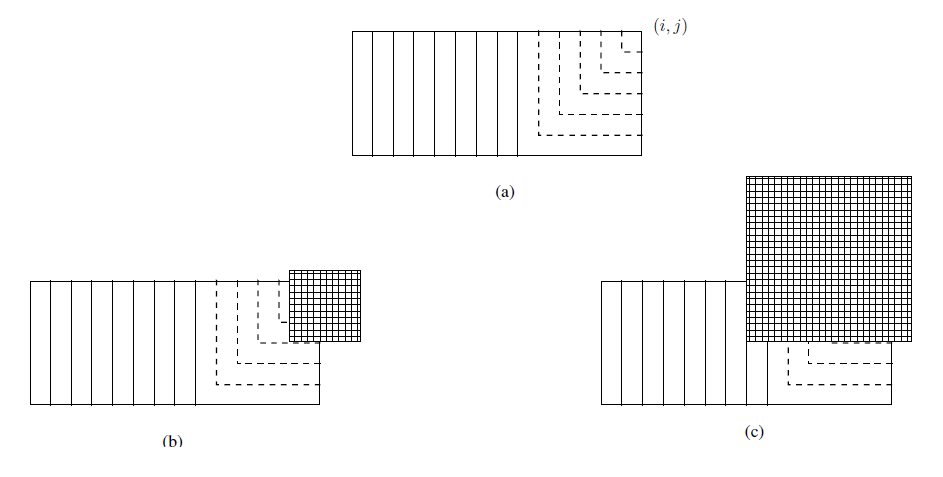}
		\caption{In (a) the vertex $(i,j)$ is not covered, in (b) the vertex $(i,j)$ is  covered from only one point in the hatched region and not covered from anywhere else in the recangle, and in (c) $(i,j)$ is  covered from only one point outside the hatched region and not covered from anywhere else in the recangle.}

	\label{r:2fig}
\end{figure}

\noindent{\sc Step 2: Understanding the event $B_{i,j}$}

\vspace{.3cm}

Let $A_{i,j} := \{(i,j) \not\in {\bf C}\}$.  There are exactly three ways that $(i,j) \not \in {\bf D}$:
\begin{itemize}
\item[(i)] $(i,j) \not \in {\bf C}$, see Figure \ref{r:2fig}(a).  As noted earlier, the probability that this occurs is 
\begin{align}
\label{r:aijdef}
a_{i,j}:= \P(A_{i,j}) = \prod_{l=0}^{j-1} g_p(l)^{2l+1} 
\prod_{k=1}^{i-j} g_p(k+j-1)^{j}.
\end{align}
\item[(ii)] there exists exactly one vertex $(k,l)$ such that $(k,l)$ is open, $0 \leq t:= \max\{i - k, j-l\} \leq j-1$ and $\rho_{(k,l)} \geq t$, see Figure \ref{r:2fig}(b). The probability that this occurs is \\
$b_{i,j}:= \prod_{k=1}^{i-j}g_p(k+j-1)^{j}\sum_{t=0}^{j-1} p G(t) g_p(t)^{2t} \prod_{l\neq t, l=0}^{j-1} g_p(l)^{2l+1} $.

\item[(iii)] there exists exactly one vertex $(k,l)$ such that $(k,l)$ is open, $1 \leq k \leq i-j$, $1\leq l \leq j\}$  and $\rho_{(k,l)} \geq i-k $, see Figure \ref{r:2fig}(c). The probability that this occurs is 
\begin{align*}
c_{i,j} := &\prod_{t=0}^{j-1} g_p(t)^{2t+1}
\Big\{\sum_{k=1}^{i-j} p G(k+j-1) g_p(k+j-1)]^{j-1}\\
& \qquad \times  \prod_{l\neq k, l=1}^{i-j} g_p(k+j-1)^{j}  \Big\}.
\end{align*}
\end{itemize} 

Thus, for $i \geq j$, and from some elementary calculations,
\begin{align}
\label{PBij}
\P_p(B_{i,j}) = & a_{i,j} + b_{i,j} + c_{i,j}\nonumber \\
  = &  a_{i,j} \left(1 + p \sum_{t=0}^{i-1} \frac{G(t)}{g_p(t)}\right).  
\end{align}

\vspace{.3cm}

\noindent{\sc Step 3: Setting up the estimates}

\vspace{.3cm}

Before we proceed further we fix some quantities. 
Since $\displaystyle\liminf_{j\longrightarrow\infty} j \mathbb{P}(\rho \geq j)>0 $ we may choose $\eta>0$ such that $\displaystyle\liminf_{j\longrightarrow\infty} j \mathbb{P}(\rho \geq j)>\eta$. Also, for this $\eta$ and our fixed $p \in (0,1)$  let $a$ be such $0<e^{-p\eta}<a<1$. Now we choose $N \geq 1$ such that for all $j \geq N$ the following hold:
\begin{flalign}
\label{conditions}
(i) &\quad j \mathbb P (\rho \geq j) > \eta,&&\nonumber \\
(ii) & \quad (1-p\eta j^{-1})^{j}<a,&&\nonumber \\
(iii) & \quad pj \eta > 1.&&
\end{flalign}
Note that (ii) above guarantees that for all $j \geq N$, we have
$g_p(j)^j < a$.

From the proof of Proposition 3.2(b) of Athreya, Roy and Sarkar (2004) we know that, for   $j \geq N$,
$$
\sum_{i=1}^\infty a_{i,j} = \sum_{i=1}^\infty \mathbb{P}_p(A_{i,j}) < \infty.
$$
Thus, from (\ref{PBij}) if we show that, for  $j \geq N$,
\begin{align}
\label{ratio11}
\sum_{i=1}^\infty  a_{i,j}\sum_{t=0}^{i-1} \frac{G(t)}{g_p(t)} < \infty
\end{align}
then we will have, for fixed $j$,  
\begin{align}
\label{Bijfinite}
\sum_{i=1}^\infty \mathbb{P}_p(B_{i,j}) < \infty \text{ whenever } j \geq N.
\end{align}

\vspace{.3cm}

\noindent{\sc  Step 4: The sum in (\ref{ratio11})}

\vspace{.3cm}

Note that, $\displaystyle\sum_{t=0}^{i-1} \frac{G(t)}{g_p(t)}=\displaystyle\sum_{t=0}^{j-1} \frac{G(t)}{g_p(t)}+ \displaystyle\sum_{k=1}^{i-j} \frac{G(k+j-1)}{g_p(k+j-1)}$. For  fixed $j \geq N$,
\begin{align}
\label{Nname}
\lefteqn{\sum_{i=j}^\infty a_{i,j}\sum_{k=1}^{i-j} \frac{G(k+j-1)}{g_p(k+j-1)}}\nonumber\\
&=\sum_{i=j}^{\infty}\prod_{t=0}^{j-1} g_p(t)^{2t+1}
\prod_{l=1}^{i-j} g_p(l+j-1)^{j}\sum_{k=1}^{i-j} \frac{G(k+j-1)}{g_(k+j-1)}\nonumber\\
&=\prod_{t=0}^{j-1} g_p(t)^{2t+1}\sum_{i=j}^{\infty} \prod_{l=j}^{i-1} g_p(l)^{j-1} \left\{\sum_{k=j}^{i-1} G(k) \prod_{h\neq k,h=j}^{i-1} g_p(h)   \right\}\nonumber\\
&=\prod_{t=0}^{j-1} g_p(t)^{2t+1}\sum_{k=j}^{\infty} \sum_{i=k+1}^{\infty}\prod_{l=j}^{i-1} g_p(l)^{j-1}  G(k) \prod_{h\neq k,h=j}^{i-1} g_p(h).   
\end{align}
Taking
\begin{align*}
e_k:=& \sum_{i=k+1}^{\infty}\prod_{l=j}^{i-1} g_p(l)^{j-1}  G(k) \prod_{h\neq k,h=j}^{i-1} g_p(h)
\end{align*}
as the inner sum in (\ref{Nname}), and breaking the sum as $i=k+1$ and $i \geq k+2$, we have
\begin{align*}
e_k =& G(k) \prod_{l=j}^{k} g_p(l)^{j-1}  \prod_{h=j}^{k-1} g_p(h) \Big( 1 + \sum_{i=k+2}^{\infty}
\prod_{l=k+1}^{i-1} g_p(l)^{j-1}  \prod_{h=k+1}^{i-1} g_p(h)\Big)\nonumber\\
& =G(k)\prod_{l=j}^{k} g_p(l)^{j-1}  \prod_{h=j}^{k-1} g_p(h)\nonumber\\
& \qquad \times \left(1+ g_p(k+1)^{j}  
 \left[1+\sum_{i=k+3}^{\infty}\prod_{l=k+2}^{i-1} g_p(l)^{j-1}  \prod_{h=k+2}^{i-1} g_p(h)\right]\right)\nonumber\\
& = G(k)\prod_{l=j}^{k} g_p(l)^{j-1}  \prod_{h=j}^{k-1} g_p(h)
\left( 1+ C(k,j) g_p(k+1)^{j}\right) ,
\end{align*}
where 
$$
C(k,j) : = 1+\sum_{i=k+3}^{\infty}\prod_{l=k+2}^{i-1} g_p(l)^{j-1} 
 \prod_{h=k+2}^{i-1} g_p(h).
 $$
Similarly, for $e_{k+1}$, as in the term in the first equality above, we have
\begin{align*}
e_{k+1} = 
 &G(k+1)\prod_{l=j}^{k+1} g_p(l)^{j-1}  \prod_{h=j}^{k} g_p(h)\nonumber\\
& \mbox{ $\;$} \qquad\times \Bigg(1+\sum_{i=k+3}^{\infty}\prod_{l=k+2}^{i-1} g_p(l)^{j-1}  \prod_{h=k+2}^{i-1} g_p(h)\Bigg)\nonumber\\
= & G(k+1)\prod_{l=j}^{k+1} g_p(l)^{j-1}  \prod_{h=j}^{k} g_p(h) C(k,j) .
\end{align*}
We see that
 \begin{align}
\label{E_kratio}
\frac{e_{k+1}}{e_k} & =C(k,j)\frac{G(k+1)}{G(k)}
\frac{g_p(k+1)^{j-1}g_p(k)}
{ 1+ C(k,j)g_p(k+1)^{j}  }\nonumber\\
& \leq \frac{C(k,j)}{ 1+ C(k,j)g_p(k+1)^{j}  }
\end{align}
because $0 \leq \frac{G(k+1)}{G(k)} \leq 1$ and $0 < g_p(k) \leq 1$ for all $k \geq 0$.
Also,  for fixed $j$, 
\begin{itemize}
\item[(i)] for all $k$ large enough,  $0.9 < g_p(k+1)^{j}\leq 1$;
\item[(ii)] for $j \geq N$ (as chosen for (\ref{conditions})), from equation (3.5) of  Athreya, Roy and Sarkar (2004), we have 
$\sum_{i=k+3}^{\infty}\prod_{l=k+2}^{i-1} g_p(l)^{j-1}  < \infty$
and hence 
$$\displaystyle\sum_{i=k+3}^{\infty}\displaystyle\prod_{l=k+2}^{i-1} g_p(l)^{j-1} \displaystyle\prod_{h=k+2}^{i-1} g_p(h) <\infty;$$
which ensures that,  for all $k$ large enough, 
$ 1 \leq C(k,j) < 1.1$.
\end{itemize}
Thus, for $j \geq N$ and all $k$ large enough, we have
$
\frac{e_{k+1}}{e_k}  < \frac{1.1}{1.9},
$
and so, by ratio test, $\sum_{k=1}^\infty e_k < \infty$. This shows, from (\ref{Nname}), that (\ref{ratio11}) and thereby (\ref{Bijfinite}) hold.

\vspace{.3cm}

\noindent{\sc Step 5: Understanding $\sum_{i,j \geq N}\mathbb{P}_p(B_{i,j})$}

\vspace{.3cm}

Now we show that, for $N$ as above,  $\sum_{i,j \geq N}\mathbb{P}_p(B_{i,j}) < \infty$. Towards this, we first observe that, for $i,j \geq 1$, by symmetry we have $\mathbb{P}_p(B_{i,j}) =  \mathbb{P}_p(B_{j,i})$, thus we need to show
\begin{align}
\label{r:doublesum}
\sum_{i,j \geq N}\mathbb{P}_p(B_{i,j})  = 2 \sum_{i=N+1}^\infty \sum_{j=N}^{i-1}\mathbb{P}_p(B_{i,j}) +
\sum_{i=N}^\infty \mathbb{P}_p(B_{i,i}) 
 < \infty.
\end{align}We will show separately that
\begin{align*}
 \sum_{i=N+1}^\infty \sum_{j=N}^{i-1}\mathbb{P}_p(B_{i,j})   < \infty  \mbox{ and }
\sum_{i=N}^\infty \mathbb{P}_p(B_{i,i})
 < \infty.
\end{align*}

Noting that $i  > j$ in the first the sum above and $i=j$ in the next sum,  we have, from the argument leading to (\ref{PBij}), 
\begin{align*}
\sum_{i=N+1}^\infty \sum_{j=N}^{i-1}\mathbb{P}_p(B_{i,j})& =\sum_{i=N+1}^\infty \sum_{j=N}^{i-1}\left[1+p\sum_{t=0}^{i-1} \frac{G(t)}{g_p(t)}\right]\mathbb{P}_p(A_{i,j}),\\
\sum_{i=N}^\infty \mathbb{P}_p(B_{i,i}) & = \sum_{i=N}^\infty  \left[1+p\sum_{t=0}^{i-1} \frac{G(t)}{g_p(t)}\right]\mathbb{P}_p(A_{i,i}). 
\end{align*}

From the proof of Proposition 3.2(b) of  Athreya, Roy and Sarkar (2004) we know that, 
$$
\sum_{i=N+1}^\infty \sum_{j=N}^{i-1}\mathbb{P}_p(A_{i,j}) < \infty \mbox{ and } \sum_{i=N}^\infty\mathbb{P}_p(A_{i,i}) < \infty, \text{ for all } p > 0;
$$
thus we need to show that
\begin{align}
\label{doublesum2}
 \sum_{i=N+1}^\infty \sum_{j=N}^{i-1}\sum_{t=0}^{i-1} \frac{G(t)}{g_p(t)}\mathbb{P}_p(A_{i,j}) < \infty
  \mbox{ and }
\sum_{i=N}^\infty \sum_{t=0}^{i-1} \frac{G(t)}{g_p(t)}\mathbb{P}_p(A_{i,i}) 
 < \infty.
\end{align}

\vspace{.3cm}

\noindent{\sc Step 6: The first sum in (\ref{doublesum2})}

\vspace{.3cm}

For the first sum, interchanging the order of the summations we have
$$
 \sum_{i=N+1}^\infty \sum_{j=N}^{i-1}\sum_{t=0}^{i-1} \frac{G(t)}{g_p(t)}\mathbb{P}_p(A_{i,j})   = \sum_{t=0}^{\infty}\frac{G(t)}{g_p(t)}\sum_{j = N}^{\infty}\sum_{i=\max(t+1,j+1)}^{\infty}\mathbb{P}_{p}(A_{i,j}) .
$$
Breaking up the inner sum according to the values  taken by $t$ in the expression on the right side above, we have
\begin{align}
\label{f:1}
&\sum_{t=0}^{\infty}\frac{G(t)}{g_p(t)}\sum_{j = N}^{\infty}\sum_{i=\max(t+1,j+1)}^{\infty}\mathbb{P}_{p}(A_{i,j})\nonumber \\
&=\sum_{t=0}^{N-1}\frac{G(t)}{g_p(t)}\sum_{j = N}^{\infty}\sum_{i=j+1}^{\infty}\mathbb{P}_{p}(A_{i,j})+\sum_{t=N}^{\infty}\frac{G(t)}{g_p(t)}\sum_{j = N}^{\infty}\sum_{i=\max(t+1,j+1)}^{\infty}\mathbb{P}_{p}(A_{i,j}).
\end{align}
From the proof of Proposition 3.2(b) of  Athreya, Roy and Sarkar (2004) we know that, 
$$
\sum_{j=N}^\infty \sum_{i=j+1}^{i-1}\mathbb{P}_p(A_{i,j}) < \infty ,$$
and hence the first term in the right side of (\ref{f:1}) is finite. 
Therefore, to show the first part of (\ref{doublesum2}),  we need to show that 
$$\sum_{t=N}^{\infty}\frac{G(t)}{g_p(t)}\sum_{j = N}^{\infty}\sum_{i=\max(t+1,j+1)}^{\infty}\mathbb{P}_{p}(A_{i,j}) < \infty.$$
Now $\sum_{j = N}^{\infty}\sum_{i=\max(t+1,j+1)}^{\infty}\mathbb{P}_{p}(A_{i,j})
= \sum_{m = 0}^{\infty}\sum_{r=\max(t-N+1,m+1)}^{\infty}\mathbb{P}_{p}(A_{N+r,N+m})$,  and breaking up the summation according to the values taken by $m$, we have
\begin{align*}
&\sum_{j = N}^{\infty}\sum_{i=\max(t+1,j+1)}^{\infty}\mathbb{P}_{p}(A_{i,j})\\
& = \sum_{m = 0}^{t-N}\sum_{r=t-N+1}^{\infty}\mathbb{P}_{p}(A_{N+r,N+m})
+ \sum_{m=t-N+1}^{\infty}\sum_{r=m+1}^{\infty}\mathbb{P}_{p}(A_{N+r,N+m})\\
= & \sum_{m=0}^{t-N}\sum_{r=t+1-N}^{\infty}\prod_{l=0}^{N+m-1}g_p(l)^{2l+1} \prod_{k=1}^{r-m}g_p(k+N+m-1)^{N+m}\\
&+\sum_{m=t-N+1}^{\infty}\sum_{r=m+1}^{\infty}\prod_{l=0}^{N+m-1}g_p(l)^{2l+1} \prod_{k=1}^{r-m}g_p(k+N+m-1)^{N+m},
\end{align*}
where we have used the expression for $\P(A_{i,j})$ as given in (\ref{r:aijdef}).

To simplify the expressions we take $\sigma_t:=\frac{G(t)}{g_p(t)}$ and
$s_m:=\displaystyle\prod_{l=0}^{N+m-1}g_p(l)^{2l+1}$. Using this notation, from the previous two equations we have
\begin{align}
\label{twosums}
& \sum_{t=N}^{\infty}\frac{G(t)}{g_p(t)}\sum_{j = N}^{\infty}\sum_{i=\max(t+1,j+1)}^{\infty}\mathbb{P}_{p}(A_{i,j}) \nonumber\\
 &\quad =   \sum_{t=N}^{\infty}\sigma_t\sum_{m=0}^{t-N}\sum_{r=t+1-N}^{\infty}s_m\prod_{k=1}^{r-m}g_p(k+N+m-1)^{N+m} \nonumber\\
 & \qquad+ \sum_{t=N}^{\infty}\sigma_t\sum_{m=t-N+1}^{\infty}\sum_{r=m+1}^{\infty}s_m\prod_{k=1}^{r-m}g_p(k+N+m-1)^{N+m}.
\end{align}
We start with the first term on the right in the above equation. Reordering the sums, we have
\begin{align}
\label{firstsum1}
&\sum_{t=N}^{\infty}\sigma_t\sum_{m=0}^{t-N}\sum_{r=t+1-N}^{\infty}s_m\prod_{k=1}^{r-m}g_p(k+N+m-1)^{N+m}\nonumber\\
& \quad = \sum_{m=0}^{\infty}s_m\sum_{t=m+N}^{\infty}\sigma_t\sum_{r=t+1-N}^{\infty}\prod_{k=1}^{r-m}g_p(k+N+m-1)^{N+m}.
\end{align}

Let 
\[\alpha_{m}:=s_m\sum_{t=m+N}^{\infty}\sigma_t\sum_{r=t+1-N}^{\infty}\prod_{k=1}^{r-m}g_p(k+N+m-1)^{N+m}\]
denote the summand.
Observe that
\begin{align}
\label{r:numdem}
&\frac{\alpha_{m+1}}{\alpha_{m}}=\frac{s_{m+1}}{s_m} \frac{\displaystyle\sum_{t=m+1+N}^{\infty}\sigma_t\displaystyle\sum_{r=t+1-N}^{\infty}\displaystyle\prod_{k=1}^{r-m-1}g_p(k+N+m)^{N+m+1}}{\displaystyle\sum_{t=m+N}^{\infty}\sigma_t\displaystyle\sum_{r=t+1-N}^{\infty}\displaystyle\prod_{k=1}^{r-m}g_p(k+N+m-1)^{N+m}}\nonumber\\
&=g_p(N+m)^{2(N+m)+1}\frac{\displaystyle\sum_{t=m+1+N}^{\infty}\sigma_t\displaystyle\sum_{r=t+1-N}^{\infty}\displaystyle\prod_{k=1}^{r-m-1}g_p(k+N+m)^{N+m+1}}{\displaystyle\sum_{t=m+N}^{\infty}\sigma_t\displaystyle\sum_{r=t+1-N}^{\infty}\displaystyle\prod_{k=1}^{r-m}g_p(k+N+m-1)^{N+m}}\nonumber\\
&= \frac{g_p(N+m)^{N+m}\displaystyle\sum_{t=m+1+N}^{\infty}\sigma_t\displaystyle\sum_{r=t+1-N}^{\infty}\displaystyle\prod_{k=1}^{r-m-1}g_p(k+N+m)^{N+m+1}}{g_p(N+m)^{-(N+m+1)}\displaystyle\sum_{t=m+N}^{\infty}\sigma_t\displaystyle\sum_{r=t+1-N}^{\infty}\displaystyle\prod_{k=1}^{r-m}g_p(k+N+m-1)^{N+m}}\, .
\end{align}
For the denominator above, we see
\begin{align*}
&g_p(N+m)^{-(N+m+1)}\displaystyle\sum_{t=m+N}^{\infty}\sigma_t\displaystyle\sum_{r=t+1-N}^{\infty}\displaystyle\prod_{k=1}^{r-m}g_p(k+N+m-1)^{N+m} \\
&=  \sigma_{m+N}g_p(N+m)^{-(N+m+1)}\displaystyle\sum_{r=m+1}^{\infty}\displaystyle\prod_{k=1}^{r-m}g_p(k+N+m-1)^{N+m}
\\
 &
+\displaystyle\sum_{t=m+1+N}^{\infty}\sigma_t\displaystyle\sum_{r=t+1-N}^{\infty}\displaystyle\prod_{k=1}^{r-m-1}g_p(k+N+m)^{N+m+1}\\
> & \displaystyle\sum_{t=m+1+N}^{\infty}\sigma_t\displaystyle\sum_{r=t+1-N}^{\infty}\displaystyle\prod_{k=1}^{r-m-1}g_p(k+N+m)^{N+m+1}.
\end{align*}

While, for the numerator of (\ref{r:numdem}), noting that  our choice of $N$ and $0 < a < 1$ as in (\ref{conditions}), implies 
$g_p(N+m)^{N+m} < a$, we have
\begin{align*}
& g_p(N+m)^{N+m}\displaystyle\sum_{t=m+1+N}^{\infty}\sigma_t\displaystyle\sum_{r=t+1-N}^{\infty}\displaystyle\prod_{k=1}^{r-m-1}g_p(k+N+m)^{N+m+1}\\
  & \qquad < \quad  a \displaystyle\sum_{t=m+1+N}^{\infty}\sigma_t\displaystyle\sum_{r=t+1-N}^{\infty}\displaystyle\prod_{k=1}^{r-m-1}g_p(k+N+m)^{N+m+1},
\end{align*}
with the strict inequality above holding if
\begin{align}
\label{r:sumbul}
\sum_{t=m+1+N}^{\infty}\sigma_t\sum_{r=t+1-N}^{\infty}\prod_{k=1}^{r-m-1}g_p(k+N+m)^{N+m+1} < \infty.
\end{align}
Thus, if (\ref{r:sumbul}) holds,  then, from (\ref{r:numdem}), we have $\frac{\alpha_{m+1}}{\alpha_{m}} < a$; and so an application of the ratio test yields that the sum in (\ref{firstsum1}) is finite.

To show (\ref{r:sumbul}) we again apply a ratio  test. First we recall that from the proof of Proposition 3.2(b) of  Athreya, Roy and Sarkar (2004) we have
\begin{align}
 \label{r:ARS3.2}
 \sum_{i=N}^\infty \prod_{k=1}^{i-j} g_p(k+j)^j < \infty \mbox{ for } N \mbox{ as in our choice.}
\end{align}
Let $\tau_t : = \sigma_t\sum_{r=t+1-N}^{\infty}\prod_{k=1}^{r-m-1}g_p(k+N+m)^{N+m+1}$. 
From (\ref{r:ARS3.2}) we have $\tau_t < \infty$.

Also
\begin{align*}
& \frac{\tau_{t+1}}{\tau_t} = \frac{G(t+1)}{G(t)}g_p(t)\\
& \qquad \quad\times \frac{g_p(t+1)^{N+m}\left(1+\displaystyle\sum_{r=t+2-N-m}^{\infty}\displaystyle\prod_{k=t+2-N-m}^{r}g_p(k+N+m)^{N+m+1}\right)}{1+g_p(t+1)^{N+m+1}\left(1+\displaystyle\sum_{r=t+2-N-m}^{\infty}\displaystyle\prod_{k=t+2-N-m}^{r}g_p(k+N+m)^{N+m+1}\right)}.
\end{align*}
From (\ref{r:ARS3.2}) we have $\displaystyle\sum_{r=t+2-N-m}^{\infty}\displaystyle\prod_{k=t+2-N-m}^{r}g_p(k+N+m)^{N+m+1} <\infty$  and so we may obtain a $t_0$ such that, for all $t \geq t_0$, 
\begin{itemize}
\item[(i)] $1 \leq 1+\displaystyle\sum_{r=t+2-N-m}^{\infty}\displaystyle\prod_{k=t+2-N-m}^{r}g_p(k+N+m)^{N+m+1} < 1.1$ and 
\item[(ii)] $g_p(t+1)^{N+m+1}\left(1+\displaystyle\sum_{r=t+2-N-m}^{\infty}\displaystyle\prod_{k=t+2-N-m}^{r}g_p(k+N+m)^{N+m+1}\right) > 0.9$. 
\end{itemize}
This choice of $t_0$ ensures that for all $t \geq t_0$, (\ref{r:sumbul}) holds.
Thus the first term on the right of (\ref{twosums}) is finite.

A similar calculation and a use of ratio test shows that the second term on the right of (\ref{twosums}) is finite. This shows that the first sum in (\ref{doublesum2}) is finite.

\vspace{.3cm}

\noindent{\sc Step 7: The second sum in (\ref{doublesum2})}

\vspace{.3cm}

Reordering the second sum in (\ref{doublesum2}) and using the notation we introduced earlier, we have
\begin{align}
\label{ssum}
\sum_{i=N}^\infty \sum_{t=0}^{i-1} \frac{G(t)}{g_p(t)}\mathbb{P}_p(A_{i,i}) = \sum_{t=0}^{N-1} \sum_{i=N}^\infty \sigma_t \mathbb{P}_p(A_{i,i}) + \sum_{t=N}^{\infty} \sum_{i=t+1}^\infty \sigma_t \mathbb{P}_p(A_{i,i}).
\end{align}
From the end of Section 3.1 of Athreya, Roy and Sarkar (2004) we know that with $N$ as above, $\sum_{i=N}^\infty \sigma_t \mathbb{P}_p(A_{i,i}) < \infty$ and hence  $\sum_{t=0}^{N-1} \sum_{i=N}^\infty \sigma_t \mathbb{P}_p(A_{i,i}) < \infty$.

Expanding the term $\mathbb{P}_p(A_{i,i})$, we have
\begin{align}
\label{mine}
\sum_{t=N}^{\infty} \sum_{i=t+1}^\infty \sigma_t \mathbb{P}_p(A_{i,i}) &= \sum_{t=N}^{\infty} \sum_{i=t+1}^\infty \sigma_t \prod_{l=0}^{i-1}g_p(l)^{2l+1}
\end{align}
Taking $a_t := \sigma_t  \sum_{i=t+1}^\infty\prod_{l=0}^{i-1}g_p(l)^{2l+1}$ and $l_t:= \sum_{r=t+2}^{\infty}\prod_{l=t+2}^{r}g_p(l)^{2l+1}$ , we see that
$$
a_t := \sigma_t\prod_{l=1}^{t}g_p(l)^{2l+1}\left\{1+g_p(t+1)^{2t+3}(1+l_t)\right\};
$$
from which we have
\begin{align*}
\frac{a_{t+1}}{a_t}& = \frac{G(t+1)}{G(t)} \frac{ g_p(t)g_p(t+1)^{2t+3}(1 + l_t)}{g_p(t+1)\left(1 +g_p(t+1)^{2t+3}(1+l_t)\right)}\\
& = \frac{G(t+1)}{G(t)} \frac{ g_p(t)g_p(t+1)^{2t+2}(1 + l_t)}{1 +g_p(t+1)^{2t+3}(1+l_t)}\\
& \leq \frac{G(t+1)}{G(t)} \frac{ g_p(t+1)^{2t+2}(1 + l_t)}{1 +g_p(t+1)^{2t+3}(1+l_t)}
\end{align*}
Also,  with $a$ as in (\ref{conditions}), 
\begin{align*}
\frac{\prod_{l=t+2}^{r+1}g_p(l)^{2l+1}}{\prod_{l=t+2}^{r}g_p(l)^{2l+1}}
& = g_p(r+1)^{2r+3}\\
&<\left[1-p\frac{p\eta}{r+1}\right]^{2r+3}\\
&<a  \text{ for } r \geq N,\\
\end{align*}
so, by the ratio test $l_t < \infty$ for $t \geq N$.
Now, since $l_t \to 0$ as $t \to \infty$, from the remark following (\ref{conditions}) we may obtain a $t_0$ such that for $t \geq t_0$, we have $g_p(t+1)^{2t+2}(1 + l_t) < a < 1$. Also  $\frac{G(t+1)}{G(t)} < 1$ for all $t$, thus 
$$ 
\frac{a_{t+1}}{a_t} < a \text{ for all } t \geq t_0,
$$
and by ratio test we have that the sum in (\ref{ssum}) is finite.

This completes the proof of Proposition \ref{r:even_main2} (ii). \hfill $\square$

\vspace{.3cm}

\noindent{{\bf Proof of Proposition \ref{r:even_main22}:}} Proposition \ref{r:even_main11} studied the 1-dimensional case of this proposition. Since the proof is along the same lines for $d \geq 2$, we provide a sketch of the argument for $d=2$.

Let $L_n:= \{(k,l) : \text{ either (i) } k=n \text{ and } 1 \leq l \leq n \text{ or (ii) } l=n \text{ and } 1 \leq k \leq n\}$ and $A_n = \{ - {\mathbf 1}, {\mathbf 0} \in ({\mathbf i}+[-\rho_{\mathbf i}, 0]^2) \text{ for some } {\mathbf i}\in L_n \text{ with } X_{\mathbf i}=1\}$.
Clearly
$
\P_p(A_n) = 1-\left(1-p\P(\rho \geq n+1)\right)^{2n-1} \sim p(2n-1) \P(\rho \geq n+1)$, so that $\sum_n \P(A_n) = \infty$ if and only if $\sum_n (2n-1) \P(\rho \geq n+1) = \infty$, i.e. if and only if $\E(\rho^2) = \infty$. Thus, by an application of the Borel-Cantelli lemma, we have  $\P_p ({\mathbf D}_{\text{rev}} = \N^2) = 1$ whenever $\E( \rho^2) = \infty$.

Conversely, for $t \in \N$,  let $\bar{{\mathbf t}}= (t, t)$ and $B_m := \{ \bar{{\mathbf t}}\in {\mathbf C}_{\rm{rev}} \text{ for all } t \geq m\}$. We will show that, if $\E(\rho^2) <  \infty$, the event $B:= \{\text{there exists } M \geq 1\text{ such that } B_m \text{ occurs for all } m \geq M\}$ is a null event, thereby implying that 
$$\P \{ \text{there exist } t_1, t_2, \ldots \text{ with } t_i \uparrow \infty \text{ such that } \bar{{\mathbf t}}_k \not \in {\mathbf C}_{\text{rev}}\} = 1$$
and hence the firework process does not percolate.

Towards this we first note that $B$ is a tail event vis-\`{a}-vis the collections of random variables $\{\chi_k\}_{k \geq 1}$ where $ \chi_k = \{X_{\mathbf i} : \mathbf i  = (i_1, i_2)\in \N^2 \text{ with } i_1, i_2 \geq k \}$. Moreover, for any
$\bar{\mathbf t} $ as in the calculations leading to the Borel-Cantelli lemma above, we have $ \P \{ \bar{\mathbf t} \not \in {\mathbf C}_{\rm{rev}}\}  > 0$ whenever $\E(\rho^2) <  \infty$. This completes the proof of Proposition \ref{r:even_main22}. \hfill 
 $\square$

\section{Poisson Boolean model}
A comparison argument with the discrete process, as in Athreya, Roy and Sarkar (2004) gives the us the following results for the Poisson Boolean model.

Let $x_1, x_2, \ldots$ be a Poisson process on $(0,\infty)^d$ of intensity $\lambda$ and $\rho, \rho_1, \rho_2, \ldots$ be i.i.d. random variables taking values in $(0,\infty)$. Let
\begin{align*}
D:= &\{z \in (0,\infty)^d: \text{there exist distinct } x_i \text{ and } x_j \text{ with } \\
& \qquad z \in (x_i + [0,\rho_i]^d ) \cap (x_j + [0,\rho_j]^d)\}.
\end{align*}
We say that $\mathbb R_+^d$ is eventually doubly covered if there exists ${\bf t} \in \mathbb R_+^d$ such that 
${\bf t} + \mathbb R_+^d \subseteq D$.
\begin{proposition}
\label{poisson1} Let $(X, \lambda, \rho)$ be a Poisson Boolean model on $\mathbb R_+^d$ with 
$$
 l:=\liminf_{x\to \infty} x \mathbb P(\rho > x)\text{  and }\; L:=\limsup_{x\to \infty} x \mathbb P(\rho > x). 
$$
(i) For $d=1$, suppose $0 < l \leq L < \infty$. There exists $1/L\leq \lambda_0 \leq 1/l$ such that
\begin{align*}
\mathbb P(\mathbb R_+ \text{ is eventually doubly covered}) &= \begin{cases}
0 & \text{if } \lambda < \lambda_0\\
1 & \text{if } \lambda > \lambda_0.
\end{cases}
\end{align*}

(ii) For $d \geq 2$
\begin{align*}
\mathbb P(\mathbb R_+^d \text{ is eventually doubly covered}) &= \begin{cases}
0 & \text{if } L =0\\
1 & \text{if } l > 0.
\end{cases}
\end{align*}
\end{proposition}

\begin{acknowledgements}
The authors acknowledge the suggestions of the referees which led to a considerable improvement of the paper. Rahul Roy also acknowledges the grant MTR/2017/000141 from DST which supported this research.
\end{acknowledgements}

 \noindent{\bf References}
 
  \vspace{.4cm}
  
 \noindent 1. {\sc Athreya, S.,  Roy, R. and Sarkar, A.} (2004).  On the coverage of space by random sets. {\it Adv. Appl. Probab. (SGSA)}\/, {\bf 36}, 1--18.\\
 
 \noindent 2. {\sc Bertacchi, D. and Zucca, F.} (2013). Rumor processes in random environment on $\N$
and on Galton-Watson trees. {\it J. Stat. Phys.}\/, {\bf 153},  486--511.\\
  
 \noindent 3. {\sc Chiu, S.N, Stoyan, D., N.,  Kendall, W.S., V. and Mecke, J. } (2013). {\it  Stochastic geometry and its applications}\/ 3rd Ed., John Wiley, Chichester.\\

\noindent 4. {\sc Gilbert, E.N.} (1961). Random plane networks, {\it J. Soc. Indust. Appl.}\/, {\bf 22}, 89--103.\\

\noindent 5. {\sc Franceschetti, M. and Meester, R.} (2007). {\it Random Networks for Communication: From Statistical Physics to Information Systems}\/, Cambridge University Press, Cambridge.\\

\noindent 6. {\sc Gallo, S., Garcia, N., Junior, V. and Rodr\'{i}guez, P.}  (2014). Rumor Processes on $\N$ and discrete
renewal processes. {\it J. Stat. Phys.}\/, {\bf 155}, 591--602.\\

\noindent 7. {\sc Gupta, P, and Kumar, P.R.} (1998).  Critical  power  for  asymptotic  connectivity  in  wireless networks. {\it Stochastic Analysis, Control, Optimization and Applications: A Volume in honor of W. H. Fleming
(W.  M.  McEneany,  G.  Yin and  Q.  Zhang,  eds.)}\/,  547--566. Birkh{\"a}user, Boston.\\

  \noindent 8. {\sc Hall, P.} (1988). {\it Introduction to the theory of coverage processes}\/ John Wiley, New York.\\
 
  \noindent 9. {\sc Junior, V., Machado, F. and Zuluaga, M.} (2011). Rumor processes on $\mathbb N$. {\it J. Appl. Probab.}\/, {\bf 48}, 624--636.\\

 \noindent 10. {\sc Junior, V., Machado, F. and Zuluaga, M.} (2014) The Cone Percolation on $\T^d$, {\it Brazilian Journal of Probability and Statistics}\/, {\bf 28} 367--675.\\
 
  \noindent 11. {\sc Junior, V., Machado, F. and Ravishankar, K.} (2016). The rumor percolation model and its variations. {\it 	arXiv:1612.03853}.\\

\noindent 12. {\sc  Maki, D. P.  and Thompson, M.} (1973).  {\it Mathematical Models and Applications, with Emphasis on Social, Life, and Management Sciences}\/, Prentice-Hall, Englewood Cliffs, NJ.\\

\noindent 13. {\sc Matheron, G. } (1968). Mod\`{e}le s\'{e}quential de partition al\'{e}atorie.
{\it Tech. Rep.}\/, Centre de Morpholgie Math\'{e}matique, Fontainbleau.\\

 \noindent 14. {\sc Meester, R. and Roy, R. } (1996). {\it Continuum percolation}\/. Cambridge University Press, New York.\\

 \noindent 15. {\sc Penrose, M.} (2003) {\it Random geometric graphs}\/, Oxford University Press, Oxford.\\
 
 \noindent 16. {\sc Sudbury, A.} (1985), The proportion of the population never hearing a rumour, {\it J. Appl. Probab.}\/, {\bf 22}, 443--446.

%
%

\end{document}